\documentclass[11pt]{amsart}
\usepackage{amssymb}
\usepackage{amsthm}
\usepackage{enumerate}
\usepackage[dvips]{graphics,color}

\theoremstyle{plain}
\newtheorem{Thm}{Theorem}[section]
\newtheorem{Lem}[Thm]{Lemma}
\newtheorem{Prop}[Thm]{Proposition}

\newtheorem{remark}[Thm]{Remark}

\theoremstyle{remark}

\small\normalsize
\numberwithin{equation}{section}
\def\beginpf{\noindent {\bf Proof:} \quad}
\newcommand{\boite}{\mbox{} \hfill $\square$}
\def\endpf{\boite\par\medskip}

\newcommand\CC{{\mathbb C}}

\newcommand\NN{{\mathbb N}}
\newcommand\RR{{\mathbb R}}

\newcommand\TT{{\mathbb T}}

\newcommand\HH{{\mathcal H}}

\begin{document}

\title[Integral representation]{Integral representation of the $n$-th derivative in de Branges-Rovnyak spaces and the norm convergence of its
reproducing kernel}

\author{Emmanuel Fricain, Javad Mashreghi}

\address{Universit\'e de Lyon; Universit\'e Lyon 1; Institut Camille Jordan CNRS UMR 5208; 43, boulevard du 11 Novembre 1918, F-69622 Villeurbanne } \email{fricain@math.univ-lyon1.fr}

\address{D\'epartement de math\'ematiques et de statistique,
         Universit\'e Laval,
         Qu\'ebec, QC,
         Canada G1K 7P4.}
\email{Javad.Mashreghi@mat.ulaval.ca}

\thanks{This work was supported by funds from NSERC (Canada) and the Jacques Cartier Center (France).}



\begin{abstract}
In this paper, we give an integral representation for the
boundary values of derivatives of functions of the de
Branges--Rovnyak spaces $\HH(b)$, where $b$ is in
 the unit ball of $H^\infty(\CC_+)$.
 In particular, we generalize a result of Ahern--Clark obtained for functions of the model
 spaces $K_b$, where $b$ is an inner function.
 Using hypergeometric series, we obtain a nontrivial formula of
 combinatorics for sums of binomial coefficients.
 Then we apply this formula to show the norm
 convergence of reproducing kernel $k_{\omega,n}^b$
 of the evaluation of $n$-th derivative of elements of $\HH(b)$ at the
 point $\omega$ as it tends radially to a point of the real axis.
\end{abstract}

\maketitle

\section{Introduction}
Let $\CC_+$ denote the upper half plane in the complex plane and let
$H^2(\CC_+)$ denote the usual Hardy space consisting of analytic
functions $f$ on $\CC_+$ which satisfy
\[
\|f\|_2:=\sup_{y>0}\left(\int_\RR |f(x+iy)|^2\,dx\right)^{1/2}<+\infty.
\]
P. Fatou \cite{fatou} proved that, for any function $f$ in
$H^2(\CC_+)$ and for almost all $x_0$ in $\RR$,
\[
f^*(x_0):=\lim_{t \to 0^+}f(x_0+it)
\]
exists. Moreover, we have $f^*\in L^2(\RR)$, $\mathcal F f^*=0$ on
$(-\infty,0)$, where $\mathcal F$ is the Fourier--Plancherel
transformation, and $\|f^*\|_2=\|f\|_2$. Of course the boundary
points where the radial limit exists depend on the function $f$.
However we cannot say more about the boundary behavior of a typical
element of $H^2(\CC_+)$. Then many authors, e.g. \cite{helson, ak70, ak71,Fricain-Mashreghi}, have studied this question by restricting the
class of functions. A particularly interesting class of subspaces of
$H^2(\CC_+)$ consists of de Branges--Rovnyak spaces.

For $\varphi\in L^\infty(\RR)$, let $T_\varphi$ stand for the
Toeplitz operator defined on $H^2(\CC_+)$ by
$$T_\varphi(f):=P_+(\varphi f),\qquad (f\in H^2(\CC_+)),$$
where $P_+$ denotes the orthogonal projection of $L^2(\RR)$ onto
$H^2(\CC_+)$. Then, for $\varphi\in L^\infty(\RR)$,
$\|\varphi\|_\infty\leq 1$,  the de Branges--Rovnyak space
$\HH(\varphi)$, associated with $\varphi$, consists of those
$H^2(\CC_+)$ functions which are in the range of the operator
$(Id-T_\varphi T_{\overline\varphi})^{1/2}$. It is a Hilbert space
when equipped with the inner product
\[
\langle (Id-T_\varphi T_{\overline \varphi})^{1/2}f,(Id-T_\varphi
T_{\overline \varphi})^{1/2}g\rangle_\varphi=\langle f,g
\rangle_2,
\]
where $f,g\in H^2(\CC_+)\ominus \hbox{ker }(Id-T_\varphi
T_{\overline \varphi})^{1/2}$.

These spaces (and more precisely their general vector-valued
version) appeared first in L. de Branges and J. Rovnyak
\cite{de-branges1, de-branges2} as universal model spaces for
Hilbert space contractions. As a special case, when $b$ is an inner
function (that is $|b|=1$ a.e. on $\RR$), the operator
$(Id-T_bT_{\overline b})$ is an orthogonal projection and $\HH(b)$
becomes a closed (ordinary) subspace of $H^2(\CC_+)$ which coincides
with the so-called model spaces $K_b=H^2(\CC_+)\ominus b
H^2(\CC_+)$. Thanks to the pioneer works of Sarason, e.g.
\cite{sarason}, we know that de Branges-Rovnyak spaces have an important role to be played in numerous questions of complex analysis and operator
theory. We mention a recent paper of A. Hartmann, D. Sarason and K. Seip \cite{HSS} who give a nice characterization of surjectivity of Toeplitz operator and the proof involves the de Branges-Rovnyak spaces. We also refer to works of J. Shapiro \cite{Shapiro1, Shapiro2} concerning the notion of angular derivative for holomorphic self-maps of the unit disk. See also a paper of J. Anderson and J. Rovnyak \cite{AR}, where generalized Schwarz-Pick estimates are given and 	a paper of M. Jury \cite{Jury}, where composition operators are studied by methods based on $\HH(b)$ spaces.

In the case where $b$ is an inner function, H. Helson \cite{helson}
studied the problem of analytic continuation across the boundary for
functions in $K_b$. Then, still when $b$ is an inner function, P.
Ahern and D. Clark \cite{ak70} characterized those points $x_0$ of
$\RR$ where every function $f$ of $K_b$ and all its derivatives up
to order $n$ have a radial limit. More precisely, if $b=BI_\mu$ is
the canonical factorization of the inner function $b$ into Blaschke
product $B$ associated with the sequence $(z_k)_k$ and singular inner
part $I_\mu$ associated with the singular measure $\mu$, then every
function $f\in K_b$ and its derivatives up to order $n$ have
finite radial limits at $x_0$ if and only if
\begin{eqnarray}\label{eq:interieur-ahern-clark}
\sum_k\frac{\Im{\rm{m}}(z_k)}{|x_0-z_k|^{2n+2}}+\int_\RR \,
\frac{d\mu(t)}{|t-x_0|^{2n+2}}<+\infty.
\end{eqnarray}
Recently, we \cite{Fricain-Mashreghi} gave an extension
of the preceding results of Helson and of Ahern--Clark. See also the paper of E. Fricain \cite{Fricain} where the orthogonal and Riesz basis of $\HH(b)$ spaces, which consist of reproducing kernels, are studied.

Now, using Cauchy formula, it is easy to see that if $b$ is inner, $\omega\in\CC_+$, $n$ is a non-negative integer and $f\in K_b$, then we have
\begin{eqnarray}\label{eq:fimodel}
f^{(n)}(\omega)=\int_\RR f(t) \, \overline{k_{\omega,n}^b(t)}\,dt,
\end{eqnarray}
where
\begin{eqnarray}\label{eq:noyaumodel}
\frac{k_{\omega,n}^b(z)}{n!}:=\frac i{2\pi}\displaystyle\frac {1-b(z)\displaystyle\sum_{p=0}^n\frac{\overline{b^{(p)}(\omega)}}{p!}(z-\overline\omega)^p}{(z-\overline{\omega})^{n+1}}\,,\qquad (z\in\CC_+).
\end{eqnarray}
A natural question is to ask if one can extend the formula (\ref{eq:fimodel}) at boundary points $x_0$. If $x_0$ is a real point which does not belong to the boundary spectrum of $b$, then $b$ and all functions of $K_b$ are analytic through a neighborhood of $x_0$ and then it is obvious that the formula (\ref{eq:fimodel}) is valid at the point $x_0$. On the other hand, if $x_0$ satisfies the condition (\ref{eq:interieur-ahern-clark}), then Ahern--Clark  \cite{ak70}
showed that the formula (\ref{eq:fimodel}) is still valid at the point $x_0\in\RR$. Recently, K. Dyakonov \cite{Dyak91, Dyak02} and then A. Baranov \cite{Baranov05} used this formula to get some Bernstein type inequalities in the model spaces $K_b$.

In this paper, our first goal is to obtain an analogue of formula
(\ref{eq:fimodel}) for the de Branges--Rovnyak spaces $\HH(b)$,
where $b$ is an {\it {arbitrary}} function in the unit ball of
$H^\infty(\CC_+)$ (not necessarily inner). We will provide an
integral representation for $f^{(n)}(\omega)$, $\omega \in \CC_+$,
and also show that under certain conditions the formula remains
valid if $\omega=x_0 \in \RR$. However, if one tries to generalize
techniques used in the model spaces $K_b$ in order to obtain such a
representation for the derivative of functions in $\HH(b)$, some
difficulties appear mainly due to the fact that the evaluation
functional in $\HH(b)$  (contrary to the model spaces $K_b$) is not
a usual integral operator. Nevertheless, we will overcome this
difficulty and provide an integral formula similar to
(\ref{eq:fimodel}) for functions in $\HH(b)$.

Our second goal is to prove the norm convergence of reproducing
kernels of evaluation functional of the $n$-th derivative as we
approach a boundary point. If $n=0$, for de Branges--Rovnyak spaces
of the unit disc, Sarason \cite[page 48]{sarason} showed that
\[
\|k_{z_0}^b\|_b^2=z_0\overline{b(z_0)}b'(z_0), \hspace{1cm}
(z_0\in\TT).
\]
We first obtain
\[
\|k_{x_0,n}^b\|_b^2=\frac{n!^2}{2i\pi}\sum_{p=0}^n\frac{\overline{b^{(p)}(x_0)}}{p!}\frac{b^{(2n+1-p)}(x_0)}{(2n+1-p)!},
\hspace{1cm} (x_0\in\RR),
\]
which is an analogue (and generalization) of  Sarason's formula for
the reproducing kernel of the $n$-th derivative for de
Branges--Rovnyak spaces of the upper half plane. Then we apply this
identity to show that $\|k_{\omega,n}^b-k_{x_0,n}^b\|_b\to 0$ as
$\omega$ tends radially to $x_0$. Again if $n=0$, this result is due
to Sarason. In establishing the norm convergence we naturally face
with the (nontrivial) finite sum
\begin{eqnarray}\label{eq:formule-combinatoire}
(-1)^{r+1}\sum_{p=0}^n\sum_{\ell=0}^n
(-2)^{p-\ell}\binom{r}{n-\ell}\binom{2n+1-r}{p}\binom{n-p+\ell}{\ell},
\end{eqnarray}
with $n,r\in\NN$, $0\leq r\leq 2n+1$. Using hypergeometric series we show that this sum is equal to $\pm
2^n$, where the choice of sign depends on $r$.

We mention a recent and very interesting work of V. Bolotnikov and A. Kheifets \cite{Bolotnikov} who obtained an
analogue of the classical Carath\'eodory--Julia theorem on boundary
derivatives. Using different techniques, the authors also obtained a condition which guarantees that we can write an analogue of formula (\ref{eq:fimodel}) for the de Branges--Rovnyak spaces $\HH(b)$. More precisely, this condition is
\begin{eqnarray}\label{eq:bk}
\liminf_{\omega\to x_0}\frac{\partial^{2n}}{\partial\omega^n\partial\overline{\omega}^n}\left(\frac{1-|b(\omega)|^2}{\Im\rm{m}\,\omega}\right)<+\infty
\end{eqnarray}
and it is stated that this is equivalent to the existence of the boundary Schwarz-Pick matrix at point $x_0$.  They also got the norm convergence (under their condition). Comparing condition (\ref{eq:bk})  with our condition (\ref{eq:condition-type-ahern-clark}) is under further investigation.

The plan of the paper is the following. In the next section, we give
some preliminaries concerning the de Branges-Rovnyak spaces. In the
third section, we establish some integral formulas for the $n$-th
derivatives of functions in $\HH(b)$. The fourth section contains the part of combinatorics of this paper. In particular, we show how we can compute the sum (\ref{eq:formule-combinatoire}) and get an interesting and quite surprising  formula. Finally, in the last section, we apply this formula of combinatorics to solve an important  problem of norm
convergence for the kernels $k_{\omega,n}^b$ corresponding to the
$n$-th derivative at points $\omega$ for functions in $\HH(b)$. More
precisely, we prove that $k_{\omega,n}^b$ tends in norm to
$k_{x_0,n}^b$ as $\omega$ tends radially to $x_0$.  We also get some
interesting relations between the derivatives of the function $b$ at
point $x_0$.

\section{preliminaries}

We first recall two general facts about the de Branges-Rovnyak
spaces. As a matter of fact, in \cite{sarason}, these results are formulated
for the unit disc. However, the same results with similar proofs also work for the upper half plane. The first one concerns the relation
between $\HH(b)$ and $\HH(\overline b)$. For $f\in H^2(\CC_+)$, we
have \cite[page 10]{sarason}
$$f\in\HH(b)\Longleftrightarrow T_{\overline b}f\in\HH(\overline b).$$
Moreover, if $f_1,f_2\in\HH(b)$, then
\begin{eqnarray}\label{eq:lien-Hb-Hbarb}
\langle f_1,f_2\rangle_b=\langle f_1,f_2\rangle_2+\langle
T_{\overline b}f_1,T_{\overline b}f_2\rangle_{\overline b}.
\end{eqnarray}
We also mention an integral representation for functions in
$\HH(\overline b)$ \cite[page 16]{sarason}. Let
$\rho(t):=1-|b(t)|^2$, $t\in\RR$, and let $L^2(\rho)$ stand for the
usual Hilbert space of measurable functions $f:\RR\to\CC$ with
$\|f\|_\rho<\infty$, where
$$\|f\|_\rho^2:=\int_\RR |f(t)|^2\rho(t)\,dt.$$
For each $w\in\CC_+$, the Cauchy kernel $k_w$ belongs to $L^2(\rho)$.
Hence, we define $H^2(\rho)$ to be the span in $L^2(\rho)$ of
the functions $k_w$ ($w\in\CC_+$). If $q$ is a function in
$L^2(\rho)$, then $q\rho$ is in $L^2(\RR)$, being the product of
$q\rho^{1/2}\in L^2(\RR)$ and the bounded function $\rho^{1/2}$.
Finally, we define the operator $C_\rho:L^2(\rho)\longrightarrow
H^2(\CC_+)$ by
$$C_\rho(q):=P_+(q\rho).$$
Then $C_\rho$ is a partial isometry from $L^2(\rho)$ onto
$\HH(\overline b)$ whose initial space equals to $H^2(\rho)$ and it
is an isometry if and only if $b$ is an extreme point of the unit
ball of $H^\infty(\CC_+)$.

In \cite{Fricain-Mashreghi}, we have studied the boundary behavior
of functions of the de Branges--Rovnyak  spaces and we mention some
parts of \cite[Theorem 3.1]{Fricain-Mashreghi} that we need here.

\begin{Thm}\label{thm:ahern-clark-demi-plan}
Let $b$ be in the unit ball of $H^\infty(\CC_+)$ and let
\[
b(z)=\prod_k e^{i\alpha_k}\frac{z-z_k}{z-\overline{z_k}} \,\,
\exp\left(-\frac{1}{i\pi}\int_\RR \frac
{tz+1}{(t-z)(t^2+1)}\,d\mu(t)\right) \,\,
\exp\left(\frac{1}{i\pi}\int_\RR \frac
{tz+1}{t-z}\frac{\log|b(t)|}{t^2+1}\,dt\right)
\]
 be its
canonical factorization. Then, for $x_0\in\RR$ and for a
non-negative integer $n$, the following are equivalent:
\begin{enumerate}
\item[$\mathrm{(i)}$] for every function $f\in\HH(b)$, $f(x_0+it), f'(x_0+it),\dots f^{(n)}(x_0+it)$ have finite  limits as $t\to 0^+$;
\item[$\mathrm{(ii)}$] we have
\begin{eqnarray}\label{eq:condition-type-ahern-clark}
\sum_k\frac {\Im{\rm m}(z_k)}{|x_0-z_k|^{2n+2}}+\int_\RR \frac
{d\mu(t)}{|x_0-t|^{2n+2}}+\int_\RR \frac
{|\log|b(t)||}{|x_0-t|^{2n+2}}\,dt<+\infty.
\end{eqnarray}
\end{enumerate}
\end{Thm}
For $f\in\HH(b)$, $x_0\in\RR$ and for a non-negative integer $n$, if $f^{(n)}(x_0+it)$ has a finite limit as $t\to 0^+$, then we define
\[
f^{(n)}(x_0):=\displaystyle\lim_{t\to 0^+} f^{(n)}(x_0+it).
\]
Moreover, under the condition (\ref{eq:condition-type-ahern-clark}),
we know that for $0\leq j\leq 2n+1$,
\begin{eqnarray}\label{eq:limite-radiale}
\displaystyle\lim_{t\to 0^+} b^{(j)}(x_0+it)
\end{eqnarray}
exists (see \cite[Lemma 4]{ak71}) and we denote this limit by
$b^{(j)}(x_0)$.

\begin{remark}
\rm{Let $x_0\in\RR$  and suppose that $x_0$ does not belong to the spectrum $\sigma(b)$ of $b$, which means (by definition) that, for some $\eta>0$, $b$ is analytic on $B(x_0,\eta):=\{z\in\CC:|z-x_0|<\eta\}$ and $|b(x)|=1$ on $(x_0-\eta,x_0+\eta)$. Denote by $a_p:=\frac{b^{(p)}(x_0)}{p!}$, $p\geq 0$. Since
\[
b(x)=\sum_{p=0}^\infty a_p(x-x_0)^p,\qquad x\in (x_0-\eta,x_0+\eta),
\]
we get
\[
1=|b(x)|^2=b(x)\overline{b(x)}=\sum_{r=0}^\infty c_r(x-x_0)^r,
\]
where $c_r=\displaystyle\sum_{p=0}^r a_p \, \overline{a_{r-p}}$.
Hence
\[
c_0=|a_0|^2=1\qquad\hbox{and}\qquad \sum_{p=0}^ra_p \,
\overline{a_{r-p}}=0,\quad (\forall r\geq 1).
\]
As we will see in the proof of Theorem~\ref{thm:main}, the condition (\ref{eq:condition-type-ahern-clark}) implies that
\[
|a_0|^2=1\qquad\hbox{and}\qquad \sum_{p=0}^r a_p \,
\overline{a_{r-p}}=0,\quad (1\leq r\leq n).
\]
Therefore, the condition (\ref{eq:condition-type-ahern-clark}) is
somehow a weaker version of the assumption $x_0\notin\sigma(b)$.}

\end{remark}

The next result gives a (standard) Taylor formula at a point on the boundary.
\begin{Lem}\label{Lem:Taylor}
Let $h$ be a holomorphic function in the upper-half plane $\CC_+$,
let $n$ be a non-negative integer and let $x_0\in\RR$. Assume that
$h^{(n)}$ has a radial limit at $x_0$. Then $h,h',\dots,h^{(n-1)}$
have radial limits at $x_0$ and
\[
h(\omega)=\sum_{p=0}^n\dfrac{h^{(p)}(x_0)}{p!}(\omega-x_0)^p+(\omega-x_0)^n\varepsilon(\omega),\qquad (\omega\in\CC_+),
\]
with $\displaystyle\lim_{t\to 0^+}\varepsilon(x_0+it)=0$.
\end{Lem}

\beginpf
The case $n=1$ is contained in \cite[Chap. VI]{sarason}. To
establish the general case one assumes as the induction hypothesis
that the property is true for $n-1$. Applying the induction
hypothesis to $h'$, we see that $h',h^{(2)},\dots,h^{(n)}$ have a
radial limit at $x_0$ and
\[
h'(\omega)=\sum_{p=0}^{n-1}\dfrac{h^{(p+1)}(x_0)}{p!}(\omega-x_0)^p+(\omega-x_0)^{n-1}\varepsilon_1(\omega),
\]
with $\displaystyle\lim_{t\to 0^+}\varepsilon_1(x_0+it)=0$.
Since $h'$ has a radial limit at $x_0$, by the case $n=1$, $h(x_0)=\displaystyle\lim_{t\to 0}h(x_0+it)$ exists and an application of Cauchy's theorem shows that
\[
h(\omega)=h(x_0)+\int_{[x_0,\omega]}h'(u)\,du,
\]
for all $\omega=x_0+it$, $t>0$. Hence we have
\begin{align*}
h(\omega)&=h(x_0)+\int_{[x_0,\omega]}\left(\sum_{p=0}^{n-1}\dfrac{h^{(p+1)}(x_0)}{p!}(u-x_0)^p+(u-x_0)^{n-1}\varepsilon_1(u)\right)\,du\\
&=\sum_{p=0}^n\dfrac{h^{(p)}(x_0)}{p!}(\omega-x_0)^p+\int_{[x_0,\omega]}(u-x_0)^{n-1}\varepsilon_1(u)\,du.
\end{align*}
Finally, let
\[
\varepsilon(\omega) = \frac{1}{(\omega-x_0)^n} \,\,
\int_{[x_0,\omega]}(u-x_0)^{n-1}\varepsilon_1(u)\,du.
\]
It is clear that $\displaystyle\lim_{t\to
0^+}\varepsilon(x_0+it)=0$.

\endpf

\section{Integral representations} \label{S:int-rep}
We first begin by proving an integral representation for the
derivatives of elements of $\HH(b)$ at points $\omega$ in the upper
half plane. Since $\omega$ is away from the boundary, the
representation is easy to establish. Let $b$ be a point in the unit
ball of $H^\infty(\CC_+)$. Recall that for $\omega\in\CC_+$, the
function
$$k_\omega^b(z)=\frac{i}{2\pi}\frac{1-\overline{b(\omega)}b(z)}{z-\overline \omega},\qquad (z\in\CC_+),$$
is the reproducing kernel of $\HH(b)$, that is
\begin{equation}\label{eq:noyau-reproduisant}
f(\omega)=\langle f,k_\omega^b\rangle_b,\qquad (f\in\HH(b)).
\end{equation}

Now let $\omega\in\CC_+$ and let $n$ be a non-negative integer. In
order to get an integral representation for the $n$th derivative of $f$  at point $\omega$ for functions in the de-Branges-Rovnyak spaces, we
need to introduce the following kernels
\begin{equation}\label{eq:kernel1-demiplan}
\frac{k_{\omega,n}^b(z)}{n!}:=\frac{i}{2\pi}\frac{1-b(z)\displaystyle\sum_{p=0}^n\frac{\overline{b^{(p)}(\omega)}}{p!}(z-\overline\omega)^p}{(z-\overline\omega)^{n+1}},\qquad
(z\in\CC_+),
\end{equation}
and
\begin{equation}\label{eq:kernel2-demiplan}
\frac{k_{\omega,n}^\rho(t)}{n!}:=\frac{i}{2\pi}\frac{\displaystyle\sum_{p=0}^n\frac{\overline{b^{(p)}(\omega)}}{p!}(t-\overline\omega)^p}{(t-\overline\omega)^{n+1}},\qquad
(t\in\RR).
\end{equation}
For $n=0$, we see that $k_{\omega,0}^b=k_\omega^b$ and
$k_{\omega,0}^\rho=\overline{b(\omega)}k_\omega$. Moreover, we also see that the kernel $k_{\omega,n}^b$ coincides with those of the inner case defined by formula (\ref{eq:noyaumodel}).

\begin{Prop}\label{Prop:representation-derivee-C+}
Let $b$ be a point in the unit ball of $H^\infty(\CC_+)$, let
$f\in\HH(b)$ and let $g\in H^2(\rho)$ be such that $T_{\overline
b}f=C_\rho(g)$. Then, for all $\omega\in\CC_+$ and for any non-negative
integer $n$, we have $k_{\omega,n}^b\in\HH(b)$ and
$k_{\omega,n}^\rho\in H^2(\rho)$ and
\begin{eqnarray}\label{eq:representation-derive-C+}
f^{(n)}(\omega)=\langle f,k_{\omega,n}^b\rangle_b=\int_\RR f(t)\overline{k_{\omega,n}^b(t)}\,dt
+\int_\RR g(t)\rho(t)\overline{k_{\omega,n}^\rho(t)}\,dt.
\end{eqnarray}
\end{Prop}

\beginpf
According to (\ref{eq:noyau-reproduisant}) and
(\ref{eq:lien-Hb-Hbarb}), we have
$$f(\omega)=\langle f,k_\omega^b\rangle_b=\langle f,k_\omega^b\rangle_2+\langle T_{\overline b}f,T_{\overline b}k_\omega^b\rangle_{\overline b}.$$
But using the fact that
$k_\omega^b=k_\omega-\overline{b(\omega)}bk_\omega$ and that
$T_{\overline b}k_\omega=\overline{b(\omega)}k_\omega$, we obtain
$$T_{\overline b}k_\omega^b=\overline{b(\omega)}\left(k_\omega-P_+(|b|^2k_\omega)\right)=\overline{b(\omega)}P_+\left((1-|b|^2)k_\omega\right)=\overline{b(\omega)}C_\rho(k_\omega),$$
which implies that
$$f(\omega)=\langle f,k_\omega^b\rangle_2+b(\omega)\langle C_\rho(g),C_\rho (k_\omega)\rangle_{\overline b}.$$
Since $C_\rho$ is a partial isometry from $L^2(\rho)$ onto
$\HH(\overline b)$, with initial space equals to $H^2(\rho)$, we   conclude that
$$f(\omega)=\langle f,k_\omega^b\rangle_2+b(\omega)\langle g,k_\omega\rangle_\rho=\langle f,k_{\omega,0}^b\rangle_2+\langle \rho g,k_{\omega,0}^\rho\rangle_2,$$
which gives the representation~(\ref{eq:representation-derive-C+})
for $n=0$.

Now straightforward computations show that
\[
\frac{\partial^n
k_{\omega,0}^b}{\partial{\overline\omega}^n}=k_{\omega,n}^b\quad\hbox{and}\quad
\frac{\partial^n
k_{\omega,0}^\rho}{\partial{\overline\omega}^n}=k_{\omega,n}^\rho.
\]
Since $k_{\omega,0}^b\in\HH(b)$ and $k_{\omega,0}^\rho\in
H^2(\rho)$, we thus have $k_{\omega,n}^b\in\HH(b)$ and
$k_{\omega,n}^\rho\in H^2(\rho)$, $n\geq 0$. The
representation~(\ref{eq:representation-derive-C+}) follows now by
induction and by differentiating  under the integral sign, which is
justified by the dominated convergence theorem.
\endpf

In the following, we show that (\ref{eq:representation-derive-C+})
is still valid at the boundary points  $x_0$ which satisfy
(\ref{eq:condition-type-ahern-clark}). We will need the boundary
analogues of the kernels (\ref{eq:kernel1-demiplan}) and
(\ref{eq:kernel2-demiplan}), i.e.
\begin{equation}\label{eq:kernel1-bord}
\frac{k_{x_0,n}^b(z)}{n!}:=\frac{i}{2\pi}\frac{1-b(z)\displaystyle\sum_{p=0}^n\frac{\overline{b^{(p)}(x_0)}}{p!}(z-x_0)^p}{(z-x_0)^{n+1}},\qquad
(z\in\CC_+),
\end{equation}
and
\begin{equation}\label{eq:kernel2-bord}
\frac{k_{x_0,n}^\rho(t)}{n!}:=\frac{i}{2\pi}\frac{\displaystyle\sum_{p=0}^n\frac{\overline{b^{(p)}(x_0)}}{p!}(t-x_0)^p}{(t-x_0)^{n+1}},\qquad
(t\in\RR\setminus\{x_0\}).
\end{equation}

The following result shows that, under condition (\ref{eq:condition-type-ahern-clark}), $k_{x_0,n}^b$ is the kernel function
in $\HH(b)$ for the functional of the $n$-th derivative at $x_0$.

\begin{Lem}\label{Lem1:noyau-pointde-R}
Let $b$ be a point in the unit ball of $H^\infty(\CC_+)$, let $n$ be
a non-negative integer and let $x_0\in\RR$. Assume that $x_0$
satisfies the condition (\ref{eq:condition-type-ahern-clark}). Then
$k_{x_0,n}^b\in\HH(b)$ and, for every function $f\in\HH(b)$, we have
\begin{equation}\label{eq:formuleHb}
f^{(n)}(x_0)=\langle f,k_{x_0,n}^b\rangle_b.
\end{equation}
\end{Lem}

\beginpf
According to Theorem \ref{thm:ahern-clark-demi-plan}, the condition (\ref{eq:condition-type-ahern-clark}) guarantees that, for
every function $f\in\HH(b)$, $f^{(n)}(\omega)$ tends to
$f^{(n)}(x_0)$, as $\omega$ tends radially to $x_0$. Therefore, an
application of the uniform boundedness principle shows that the
functional $f\longmapsto f^{(n)}(x_0)$ is bounded on $\HH(b)$.
Hence, by Riesz' theorem, there exists $\varphi_{x_0,n}\in\HH(b)$
such that
$$f^{(n)}(x_0)=\langle f,\varphi_{x_0,n}\rangle_b,\qquad (f\in\HH(b)).$$
Since
$$f^{(n)}(\omega)=\langle f,\frac {\partial^n k_{\omega,0}^b}{\partial{\overline\omega}^n} \rangle_b=
\langle f,k_{\omega,n}^b\rangle_b,\qquad (f\in\HH(b)),$$ we see that
$k_{\omega,n}^b$ tends weakly to $\varphi_{x_0,n}$, as $\omega$ tends
radially to $x_0$. Thus, for $z\in\CC_+$, we can write
$$\begin{aligned}
\varphi_{x_0,n}(z)&=\langle \varphi_{x_0,n},k_z^b\rangle_b=\displaystyle\lim_{t\to 0^+}\langle k_{x_0+it,n}^b,k_z^b\rangle_b=\displaystyle\lim_{t\to 0^+}  k_{x_0+it,n}^b(z)\\
&=\displaystyle\lim_{t\to 0^+}n!\frac{i}{2\pi}\frac{1-b(z)\displaystyle\sum_{p=0}^n\frac{\overline{b^{(p)}(x_0+it)}}{p!}(z-x_0+it)^p}{(z-x_0+it)^{n+1}}\\
&=n!\frac{i}{2\pi}\frac{1-b(z)\displaystyle\sum_{p=0}^n\frac{\overline{b^{(p)}(x_0)}}{p!}(z-x_0)^p}{(z-x_0)^{n+1}},
\end{aligned}$$
which implies that $\varphi_{x_0,n}=k_{x_0,n}^b$.  Hence $k_{x_0,n}^b\in\HH(b)$ and for every function $f\in\HH(b)$ we have
$$f^{(n)}(x_0)=\langle f,k_{x_0,n}^b\rangle_b.$$

\endpf

For $n=0$, Lemma~\ref{Lem1:noyau-pointde-R} appears in \cite[Chap.
V]{sarason}, in the context of the unit disc. The problem with the
representation (\ref{eq:formuleHb}) is that the inner product in
$\HH(b)$ is not an explicit integral formula and thus it is not
convenient to use it. That is why we prefer to have an integral
formula of type (\ref{eq:representation-derive-C+}).

If $x_0$ satisfies the condition
(\ref{eq:condition-type-ahern-clark}) we also have
$k_{x_0,n}^\rho\in L^2(\rho)$. Indeed, according to
(\ref{eq:kernel2-bord}), it suffices to prove that $(t-x_0)^{-j}\in
L^2(\rho)$, for $1\leq j\leq n+1$. Since $\rho\leq 1$, it is enough
to verify this fact in a neighborhood of $x_0$, say
$I_{x_0}=[x_0-1,x_0+1]$. But according to the condition
(\ref{eq:condition-type-ahern-clark}), we have
\[
\int_{I_{x_0}}\frac {1-|b(t)|^2}{|t-x_0|^{2j}}\,dt\leq
2\int_{I_{x_0}}\frac{|\log|b(t)||}{|t-x_0|^{2j}}\,dt\leq
2\int_{I_{x_0}}\frac{|\log|b(t)||}{|t-x_0|^{2(n+1)}}\,dt<+\infty.
\]

\begin{Thm}\label{thm:main}
Let $b$ be a point in the unit ball of $H^\infty(\CC_+)$, let $n$ be
a non-negative integer, let $f\in\HH(b)$ and let $g\in H^2(\rho)$ be such that $T_{\overline
b}f=C_\rho(g)$. Then, for every point $x_0\in\RR$ satisfying the condition (\ref{eq:condition-type-ahern-clark}), we have
\begin{equation}\label{eq:representation-de-Branges}
f^{(n)}(x_0)=\int_\RR f(t)\overline{k_{x_0,n}^b(t)}\,dt
+\int_\RR g(t)\rho(t)\overline{k_{x_0,n}^\rho(t)}\,dt.
\end{equation}
\end{Thm}

\beginpf
Recall that according to (\ref{eq:limite-radiale}), the condition (\ref{eq:condition-type-ahern-clark}) guarantees that $b^{(j)}(x_0)$ exists for $0\leq j\leq 2n+1$. Moreover, Lemma~\ref{Lem1:noyau-pointde-R} implies that $k_{x_0,p}^b\in\HH(b)$, for $0\leq p\leq n$. First of all, we prove that
\[
h_{x_0,n}(z):=\frac{b(z)-\displaystyle\sum_{p=0}^n\frac{b^{(p)}(x_0)}{p!}(z-x_0)^p}{(z-x_0)^{n+1}},\qquad (z\in\CC_+),
\]
satisfies
\begin{equation}\label{eq1:cle}
h_{x_0,n}=2i\pi\displaystyle\sum_{p=0}^n \frac{b^{(n-p)}(x_0)}{(n-p)!p!}k_{x_0,p}^b.
\end{equation}
To simplify a little bit the next computations, we put $a_p:=\frac{b^{(p)}(x_0)}{p!}$, $0\leq p\leq n$. According to (\ref{eq:kernel1-demiplan}), we have
\begin{align*}
2i\pi\sum_{p=0}^n a_{n-p}\frac{k_{x_0,p}^b(z)}{p!}&=\sum_{p=0}^n a_{n-p} \left(\displaystyle\frac{\displaystyle\sum_{j=0}^p \overline{a_j}(z-x_0)^jb(z)-1}{(z-x_0)^{p+1}}\right)\\
&=\frac{\displaystyle\sum_{p=0}^n a_{n-p}(z-x_0)^{n-p}\left(b(z)\displaystyle\sum_{j=0}^p \overline{a_j}(z-x_0)^j-1\right)}{(z-x_0)^{n+1}}\\
&=\frac{b(z)\left(\displaystyle\sum_{p=0}^n\sum_{j=0}^p a_{n-p}\overline{a_j}(z-x_0)^{n-p+j}\right)-\displaystyle\sum_{k=0}^n a_k(z-x_0)^k}{(z-x_0)^{n+1}}.
\end{align*}
Therefore, we see that (\ref{eq1:cle}) is equivalent to
\begin{equation}\label{eq2:cle}
\displaystyle\sum_{p=0}^n\sum_{j=0}^p a_{n-p}\overline{a_j}(z-x_0)^{n-p+j}=1.
\end{equation}
But, putting $j=\ell-n+p$, we obtain
\begin{align*}
\displaystyle\sum_{p=0}^n\sum_{j=0}^p a_{n-p}\overline{a_j}(z-x_0)^{n-p+j}&=\displaystyle\sum_{\ell=0}^n\left(\sum_{p=n-\ell}^n a_{n-p}\overline{a_{\ell-n+p}}\right)(z-x_0)^\ell\\
&=\displaystyle\sum_{\ell=0}^n\left(\displaystyle\sum_{q=0}^\ell a_{\ell-q}\overline{a_q}\right)(z-x_0)^\ell.
\end{align*}
Consequently, (\ref{eq2:cle}) is equivalent to
\begin{equation}\label{eq3:cle}
|b(x_0)|^2=1\quad\hbox{and}\quad \displaystyle\sum_{q=0}^\ell a_{\ell-q}\overline{a_q}=0,\qquad (1\leq\ell\leq n).
\end{equation}
Now if we define
\[
\varphi(z):=1-b(z)\displaystyle\sum_{p=0}^n \overline{a_p}(z-x_0)^p,\qquad (z\in\CC_+),
\]
then $\varphi$ is holomorphic in $\CC_+$ and according to (\ref{eq:limite-radiale}), $\varphi$ and its derivatives up to $2n+1$ have radial limits at $x_0$. An application of Lemma \ref{Lem:Taylor} shows that we can write
\[
\varphi(z)=\sum_{p=0}^n\dfrac{\varphi^{(p)}(x_0)}{p!}(z-x_0)^p+o((z-x_0)^n),
\]
as $z$ tends radially to $x_0$. Assume that there exists $k\in\{0,\dots,n\}$ such that $\varphi^{(k)}(x_0)\not=0$ and let
\[
j_0:=\min\{0\leq p\leq n:\varphi^{(p)}(x_0)\not=0\}.
\]
Hence, as $t\to 0^+$,
\[
|k_{x_0,n}^b(x_0+it)|\sim\frac{1}{2\pi}\frac{|\varphi^{(j_0)}(x_0)|}{j_0!}t^{j_0-(n+1)},
\]
which implies that $\displaystyle\lim_{t\to 0^+}|k_{x_0,n}^b(x_0+it)|=+\infty$. This is a contradiction with the fact that $k_{x_0,n}^b$ belongs to $\HH(b)$ and has a finite radial limit at $x_0$. Therefore we necessarily have $\varphi^{(\ell)}(x_0)=0$, $0\leq\ell\leq n$. But $\varphi(x_0)=1-b(x_0)\overline{b(x_0)}=1-|b(x_0)|^2$ and if we use the Leibniz' rule to compute the derivative of $\varphi$, for $1\leq\ell\leq n$, we get
\[
\varphi^{(\ell)}(x_0)=-\sum_{p=0}^\ell \overline{a_p}\binom{\ell}{p}p!b^{(\ell-p)}(x_0)=-\ell!\sum_{p=0}^\ell\overline{a_p}a_{\ell-p},
\]
which gives (\ref{eq3:cle}).
Hence (\ref{eq1:cle}) is proved. According to Lemma~\ref{Lem1:noyau-pointde-R},  (\ref{eq1:cle}) implies $h_{x_0,n}\in\HH(b)$. Now for almost all $t\in\RR$, we have
\begin{align*}
\overline{b(t)}\frac{k_{x_0,n}^b(t)}{n!}&=\frac{i}{2\pi}\frac{\overline{b(t)}-|b(t)|^2\displaystyle\sum_{p=0}^n \overline{a_p}(t-x_0)^p}{(t-x_0)^{n+1}}\\
&=\frac{i}{2\pi}(1-|b(t)|^2)\frac{\displaystyle\sum_{p=0}^n\overline{a_p}(t-x_0)^p}{(t-x_0)^{n+1}}+\frac{i}{2\pi}\frac{\overline{b(t)}-\displaystyle\sum_{p=0}^n \overline{a_p}(t-x_0)^p}{(t-x_0)^{n+1}}\\
&=\rho(t)\frac{k_{x_0,n}^\rho}{n!}+\frac{i}{2\pi}\overline{h_{x_0,n}(t)}.
\end{align*}
Since $h_{x_0,n}\in\HH(b)\subset H^2(\CC_+)$, we get that $P_+(\overline b k_{x_0,n}^b)=P_+(\rho k_{x_0,n}^\rho)$, which can be written as $T_{\overline b}k_{x_0,n}^b=C_\rho k_{x_0,n}^\rho$. It follows from (\ref{eq:lien-Hb-Hbarb}) and Lemma \ref{Lem1:noyau-pointde-R} that
\begin{align*}
f^{(n)}(x_0)&=\langle f,k_{x_0,n}^b\rangle_b \\
&=\langle f,k_{x_0,n}^b\rangle_2+\langle T_{\overline b}f,T_{\overline b}k_{x_0,n}^b\rangle_{\overline b}\\
&=\langle f,k_{x_0,n}^b\rangle_2+\langle g,k_{x_0,n}^\rho\rangle_{\rho}\\
&=\int_\RR f(t)\overline{k_{x_0,n}^b(t)}\,dt+\int_\RR g(t)\rho(t)\overline{k_{x_0,n}^\rho(t)}\,dt,
\end{align*}
which proves the relation (\ref{eq:representation-de-Branges}).

\endpf

If $b$ is inner, then it is clear that the second integral in
(\ref{eq:representation-de-Branges}) is zero and we obtain the
formula of Ahern--Clark  (\ref{eq:fimodel}).

\section{A formula of combinatorics}
We first recall some well-known facts concerning hypergeometric
series (see \cite{Andrews, Slater}).

The ${}_2F_1$ hypergeometric series is a power series in $z$ defined by
\begin{equation}\label{eq:gauss-definition}
{}_2F_1\!{\scriptsize{\left[\begin{matrix}a,b\\
c\end{matrix};z\right]}}=\sum_{p=0}^{+\infty}\frac{(a)_p(b)_p}{p!(c)_p}z^p,
\end{equation}
where $a,b,c \in \CC$, $c\not=0,-1,-2,\dots$,  and
\[
(t)_p:=\begin{cases}
1,& \mbox{ if } p=0,\\
t(t+1)\dots(t+p-1),&  \mbox{ if } p\geq 1.
\end{cases}
\]
We see that the hypergeometric series reduces to a polynomial of
degree $n$ in $z$ when $a$ or $b$ is equal to $-n$,
$(n=0,1,2,\dots)$. It is clear that the radius of convergence of the
${}_2F_1$ series is equal to 1. One can show  that when
$\Re\rm{e}(c-a-b)\leq -1$ this series is divergent on the entire
unit circle, when $-1<\Re\rm{e}(c-a-b)\leq 0$ this series converges
on the unit circle except for $z=1$ and when  $0<\Re\rm{e}(c-a-b)$
this series is (absolutely) convergent on the entire unit circle
(see \cite[Theorem 2.1.2]{Andrews}).

We note that a power series $\sum_p\alpha_p z^p$ $(\alpha_0=1)$ can be written as a hypergeometric series ${}_2F_1\!{\scriptsize{\left[\begin{matrix}a,b\\
c\end{matrix};z\right]}}$ if and only if
\begin{eqnarray}\label{eq:test-hypergeometrique}
\frac{\alpha_{p+1}}{\alpha_p}=\frac{(p+a)(p+b)}{(p+1)(p+c)}.
\end{eqnarray}

Finally we recall two useful well-known formulas \cite[page 68]{Andrews} for the hypergeometric  series:\begin{equation}\label{eq:euler}
{}_2F_1\!{\scriptsize{\left[\begin{matrix}a,b\\
c\end{matrix};z\right]}}=(1-z)^{c-a-b} {}_2F_1\!{\scriptsize{\left[\begin{matrix}c-a,c-b\\
c\end{matrix};z\right]}}\qquad (\hbox{Euler's formula}),
\end{equation}
and
\begin{equation}\label{eq:formule1}
{}_2F_1\!{\scriptsize{\left[\begin{matrix}a,b\\
c\end{matrix};\tfrac 12\right]}}=2^a {}_2F_1\!{\scriptsize{\left[\begin{matrix}a,c-b\\
c\end{matrix};-1\right]}},\qquad \Re{\rm{e}}(b-a)>-1,\quad (\hbox{Pfaff's formula}).
\end{equation}

Now we state the result which we use in the last section.

\begin{Prop}\label{Prop:Fred}
Let $n,r\in\NN$, $0\leq r\leq 2n+1$ and define
\begin{equation}\label{eq:definition-de-ANR}
A_{n,r}:=(-1)^{r+1}\sum_{p=0}^n\sum_{\ell=0}^n (-2)^{p-\ell}\binom{r}{n-\ell}\binom{2n+1-r}{p}\binom{n-p+\ell}{\ell}.
\end{equation}
Then
$$
A_{n,r}=\begin{cases}
-2^n,&0\leq r\leq n \\
2^n,&n+1\leq r\leq 2n+1.
\end{cases}
$$
\end{Prop}

For the proof of this result, we need the following lemma.

\begin{Lem}\label{Lem:Zeng}
For $m\in\NN$, we have
\begin{equation}\label{eq:formule2}
\sum_{k=0}^m\binom{m}{k}(z-1)^{-k}{}_2F_1\!{\scriptsize{\left[\begin{matrix}a,b-k\\
c\end{matrix};z\right]}}=\frac{(c-a)_m}{(c)_m}\left(\frac{z}{z-1}\right)^m{}_2F_1\!{\scriptsize{\left[\begin{matrix}a,b\\
c+m\end{matrix};z\right]}}.
\end{equation}
\end{Lem}

\beginpf
First note that (\ref{eq:formule2}) is equivalent to
\begin{equation}\label{eq:formule2-bis}
\sum_{k=0}^m\binom{m}{k}(1-z)^{-k}(-1)^{m-k}{}_2F_1\!{\scriptsize{\left[\begin{matrix}a,b-k\\
c\end{matrix};z\right]}}=\frac{(c-a)_m}{(c)_m}\left(\frac{z}{1-z}\right)^m{}_2F_1\!{\scriptsize{\left[\begin{matrix}a,b\\
c+m\end{matrix};z\right]}},
\end{equation}
and denote by $LH$ the left hand side of the inequality
(\ref{eq:formule2-bis}). Applying transformation (\ref{eq:euler}),
we obtain
\[
LH=\sum_{k=0}^m\binom{m}{k}(-1)^{m-k}(1-z)^{c-a-b}{}_2F_1\!{\scriptsize{\left[\begin{matrix}c-a,c-b+k\\
c\end{matrix};z\right]}}.
\]
Now we introduce the operator of difference $\Delta$ defined by
$\Delta f(x)=f(x+1)-f(x)$. Then it is well-known and easy to verify
that
\[
\Delta^mf(x)=\sum_{k=0}^m \binom{m}{k}(-1)^{m-k}f(x+k).
\]
Using this formula, we see that $LH=(1-z)^{c-a-b}\Delta^m f(c-b)$,
with
\[
f(x):={}_2F_1\!{\scriptsize{\left[\begin{matrix}c-a,x\\
c\end{matrix};z\right]}}.
\]
But now we can compute $\Delta^mf(x)$. Indeed, we have
\begin{align*}
\Delta f(x)=&\sum_{k=0}^{+\infty}\frac{(c-a)_k}{(c)_k}((x+1)_k-(x)_k)\frac{z^k}{k!}\\
=&\sum_{k=1}^{+\infty}\frac{(c-a)_k}{(c)_k}(x+1)_{k-1}\frac{z^k}{(k-1)!}\\
=&\frac{(c-a)}{c}\,z\,{}_2F_1\!{\scriptsize{\left[\begin{matrix}c-a+1,x+1\\
c+1\end{matrix};z\right]}},
\end{align*}
and by induction, it follows that
\[
\Delta^m f(x)=\frac{(c-a)_m}{(c)_m}z^m{}_2F_1\!{\scriptsize{\left[\begin{matrix}c-a+m,x+m\\
c+m\end{matrix};z\right]}}.
\]
Therefore, we get
\[
LH=(1-z)^{c-a-b}\frac{(c-a)_m}{(c)_m}z^m {}_2F_1\!{\scriptsize{\left[\begin{matrix}c-a+m,c-b+m\\
c+m\end{matrix};z\right]}}.
\]
Applying once more Euler's formula, we obtain (\ref{eq:formule2-bis}).

\endpf

{\bf{Proof of Proposition \ref{Prop:Fred}:}}$\quad$ Changing $\ell$
into $n-\ell$ in the second sum of (\ref{eq:definition-de-ANR}), we
see that
\begin{equation}\label{eq:transformation1-ANR}
A_{n,r}=(-1)^{r+1}\sum_{p=0}^n\sum_{\ell=0}^n (-2)^{p+\ell-n}\binom{r}{\ell}\binom{2n+1-r}{p}\binom{2n-p-\ell}{n-\ell}.
\end{equation}
Hence
\begin{align*}
A_{n,2n+1-r}&=(-1)^{2n+1-r+1}\sum_{p=0}^n\sum_{\ell=0}^n (-2)^{p+\ell-n}\binom{2n+1-r}{\ell}\binom{r}{p}\binom{2n-p-\ell}{n-\ell}\\
&=-(-1)^{r+1}\sum_{p=0}^n\sum_{\ell=0}^n (-2)^{p+\ell-n}\binom{2n+1-r}{\ell}\binom{r}{p}\binom{2n-p-\ell}{n-p}\\
&=-A_{n,r}.
\end{align*}
Therefore, it is sufficient to show $A_{n,r}=-2^n$ for $0\leq r\leq
n$ and then the result for $n+1\leq r\leq 2n+1$ will follow
immediately.

We will now assume that $0\leq r\leq n$. Changing $p$ to $n-p$ in the first sum of (\ref{eq:transformation1-ANR}) and permuting the two sums, we get
\begin{align*}
A_{n,r}=(-1)^{r+1}\sum_{\ell=0}^n(-2)^l\binom{r}{\ell}\sum_{p=0}^n (-2)^{-p}\binom{2n+1-r}{n-p}\binom{n+p-\ell}{n-\ell}.
\end{align*}
According to (\ref{eq:gauss-definition}) and (\ref{eq:test-hypergeometrique}), we see that
\[
\sum_{p=0}^n (-2)^{-p}\binom{2n+1-r}{n-p}\binom{n+p-\ell}{n-\ell}=\binom{2n+1-r}{n}{}_2F_1\!{\scriptsize{\left[\begin{matrix}n-\ell+1,-n\\
n+2-r\end{matrix};\tfrac 12\right]}},
\]
which implies
\[
A_{n,r}=(-1)^{r+1}\binom{2n+1-r}{n}\sum_{\ell=0}^n (-2)^\ell\binom{r}{\ell}{}_2F_1\!{\scriptsize{\left[\begin{matrix}n-\ell+1,-n\\
n+2-r\end{matrix};\tfrac 12\right]}}.
\]
Since $r\leq n$ and since $\binom{r}{\ell}=0$ if $r<\ell$, the sum (in the last equation) ends at $\ell=r$ and we can apply Lemma \ref{Lem:Zeng} with $m=r$, $a=-n$, $b=n+1$, $c=n+2-r$ and $z=\frac 12$. Therefore
\begin{align*}
A_{n,r}=&(-1)^{r+1}\binom{2n+1-r}{n}\frac{(2n+2-r)_r}{(n+2-r)_r}(-1)^r{}_2F_1\!{\scriptsize{\left[\begin{matrix}-n,n+1\\
n+2\end{matrix};\tfrac 12\right]}}\\
=&-\frac{(2n+1)!}{(n+1)!n!}\,\,{}_2F_1\!{\scriptsize{\left[\begin{matrix}-n,n+1\\
n+2\end{matrix};\tfrac 12\right]}}.
\end{align*}
We now use formula (\ref{eq:formule1}) which gives
\begin{align*}
{}_2F_1\!{\scriptsize{\left[\begin{matrix}-n,n+1\\
n+2\end{matrix};\tfrac 12\right]}}&=2^{-n}{}_2F_1\!{\scriptsize{\left[\begin{matrix}-n,1\\
n+2\end{matrix};-1\right]}}\\
&=2^{-n}\sum_{i=0}^\infty\frac{(-n)_i(1)_i}{(n+2)_i}\frac{(-1)^i}{i!}\\
&=2^{-n}\sum_{i=0}^n \binom{n}{i}\frac{i!(n+1)!}{(n+i+1)!},
\end{align*}
where we have used $(1)_i=i!$,
$\displaystyle\frac{(-n)_i(-1)^i}{i!}=\binom{n}{i}$ and
$(n+2)_i=\displaystyle\frac{(n+i+1)!}{(n+1)!}$. Hence
\[
{}_2F_1\!{\scriptsize{\left[\begin{matrix}-n,n+1\\
n+2\end{matrix};\tfrac 12\right]}}=2^{-n}(n+1)!n!\sum_{i=0}^n \frac{1}{(n+i+1)!(n-i)!}
\]
which implies
\[
A_{n,r}=-2^{-n}\sum_{i=0}^n\binom{2n+1}{i}.
\]
But
\[
\sum_{i=0}^n\binom{2n+1}{i}=\frac{1}{2}\sum_{i=0}^{2n+1}\binom{2n+1}{i}=2^{2n}
\]
which ends the proof.

\endpf

\section{Norm convergence for the reproducing kernels.}
In Section \ref{S:int-rep}, we saw that if $x_0 \in \RR$ satisfies
(\ref{eq:condition-type-ahern-clark}), then $k_{\omega,n}^b$
tends {\em weakly} to $k_{x_0,n}^b$ in $\HH(b)$ as $\omega$ approaches
radially to $x_0$. It is natural to ask if this weak convergence can
be replaced by norm convergence. In other words, is it true that
$\|k_{\omega,n}^b-k_{x_0,n}^b\|_b\to 0$ as $\omega$ tends radially
to $x_0$?

In \cite{ak70}, Ahern and Clark said that they can prove this result
for the case where $b$ is inner and $n=0$. For general functions $b$
in the unit ball of $H^\infty$, Sarason \cite[Chap. V]{sarason} got
this norm convergence for the case $n=0$. In this section, we prove
the general case.

Since we already have weak convergence, to prove the norm
convergence, it is sufficient to prove that
$\|k_{\omega,n}^b\|_b\to\|k_{x_0,n}^b\|_b$ as $\omega$ tends
radially to $x_0$. Therefore we need to compute $\|k_{x_0,n}^b\|_b$.
For $n=0$, in the context of the unit disc, Sarason \cite[Chap.
V]{sarason} proved that
$\|k_{z_0}^b\|_b^2=z_0\overline{b(z_0)}b'(z_0)$, $z_0\in\TT$. We can
give an analogue of this formula showing that the norm of
$k_{x_0,n}^b$ can be expressed in terms of the derivatives of $b$ at
$x_0$.

\begin{Prop}\label{Prop:cacul-de-norme}
Let $b$ be a point in the unit ball of $H^\infty(\CC_+)$, let $n$ be
a non-negative integer and let $x_0\in\RR$ satisfying the condition (\ref{eq:condition-type-ahern-clark}). Then
\[
\|k_{x_0,n}^b\|_b^2=\frac{n!^2}{2i\pi}\sum_{p=0}^n\frac{\overline{b^{(p)}(x_0)}}{p!}\frac{b^{(2n+1-p)}(x_0)}{(2n+1-p)!}.
\]
\end{Prop}

\beginpf Following the notations of Section 3, we define
\[
\varphi(z)=1-b(z)\sum_{p=0}^n\frac{\overline{b^{(p)}(x_0)}}{p!}(z-x_0)^p.
\]
Then, by (\ref{eq:limite-radiale}) and Lemma
\ref{Lem:Taylor}, as $z$ tends radially to $x_0$,  we have
\[
k_{x_0,n}^b(z)=\frac{in!}{2\pi}(z-x_0)^{-n-1}\left(\sum_{p=0}^{2n+1}\dfrac{\varphi^{(p)}(x_0)}{p!}(z-x_0)^p+o((z-x_0)^{2n+1})\right).
\]
As we have shown in the proof of Theorem~\ref{thm:main},
$\varphi^{(k)}(x_0)=0$ if $0\leq k\leq n$. Hence
\[
k_{x_0,n}^b(z)=\frac{in!}{2\pi}\sum_{p=0}^{n}\dfrac{\varphi^{(p+n+1)}(x_0)}{(p+n+1)!}(z-x_0)^p+o((z-x_0)^n).
\]
Using once more Lemma~\ref{Lem:Taylor}, we can also write
\[
k_{x_0,n}^b(z)=\sum_{p=0}^n\dfrac{{(k_{x_0,n}^b)}^{(p)}(x_0)}{p!}(z-x_0)^p+\,o((z-x_0)^n),
\]
which implies
\[
{(k_{x_0,n}^b)}^{(p)}(x_0)=\frac{in!}{2\pi}\frac{p!}{(p+1+n)!}\varphi^{(p+n+1)}(x_0).
\]
But, according to Lemma \ref{Lem1:noyau-pointde-R}, we have $\|k_{x_0,n}^b\|_b^2=(k_{x_0,n}^b)^{(n)}(x_0)$ and we get
\[
\|k_{x_0,n}^b\|_b^2=\frac{in!^2}{2\pi}\frac{\varphi^{(2n+1)}(x_0)}{(2n+1)!}.
\]
Finally, the result follows by Leibniz' rule.

\endpf

The next result provides an affirmative answer to the question of norm convergence.

\begin{Thm}
Let $b$ be a point in the unit ball of $H^\infty(\CC_+)$, let $n$ be
a non-negative integer and let $x_0\in\RR$ satisfying the condition (\ref{eq:condition-type-ahern-clark}). Then
\[
\left\|k_{\omega,n}^b-k_{x_0,n}^b\right\|_b\longrightarrow 0,\quad\hbox{as $\omega$ tends radially to $x_0$.}
\]
\end{Thm}

\beginpf
We denote by $a_p(\omega):=\frac{b^{(p)}(\omega)}{p!}$ and $a_p:=a_p(x_0)$. We recall that
\[
k_{\omega,n}^b(z)=\frac{in!}{2\pi}\left(
\displaystyle\frac{1}{(z-\overline\omega)^{n+1}}-\sum_{p=0}^n\overline{a_p(\omega)}(z-\overline\omega)^{p-n-1}b(z)\right).
\]
We have
\[
\frac{\partial^n}{\partial z^n}\left(\frac{1}{(z-\overline\omega)^{n+1}}\right)=(-1)^n\frac{(2n)!}{n!}\frac{1}{(z-\overline\omega)^{2n+1}},
\]
and by Leibniz' rule
\[
\frac{\partial^n}{\partial z^n}\left((z-\overline\omega)^{p-n-1}b(z)\right)=\sum_{\ell=0}^n\binom{n}{\ell}(-1)^\ell\frac{(n-p+\ell)!}{(n-p)!}(z-\overline\omega)^{p-n-\ell-1}b^{(n-\ell)}(z).
\]
According to Proposition~\ref{eq:representation-derive-C+}, we have
$\|k_{\omega,n}^b\|_b^2=(k_{\omega,n}^b)^{(n)}(\omega)$, which
implies
\begin{equation}\label{eq1:asymptot-norme}
\|k_{\omega,n}^b\|_b^2=\dfrac{in!}{2\pi}\frac{(-1)^n\dfrac{(2n)!}{n!}-\displaystyle\sum_{p=0}^n\sum_{\ell=0}^n\binom{n}{\ell}(-1)^\ell\frac{(n-p-\ell)!}{(n-p)!}(\omega-\overline\omega)^{n+p-\ell}\overline{a_p(\omega)}b^{(n-\ell)}(\omega)}{(\omega-\overline\omega)^{2n+1}}.
\end{equation}
For $0\leq s\leq n$, the function $b^{(s)}$ is analytic in the upper-half plane and its derivative of order $2n+1-s$, which coincides with $b^{(2n+1)}$, has a radial limit at $x_0$. According to Lemma \ref{Lem:Taylor}, as $\omega$ tends radially to $x_0$, we have
\[
b^{(s)}(\omega)=\sum_{r=s}^{2n+1}a_r\frac{r!}{(r-s)!}(\omega-x_0)^{r-s}+o((\omega-x_0)^{2n+1-s}).
\]
Hence if we put $\omega=x_0+it$, we get
\[
(\omega-\overline\omega)^s b^{(s)}(\omega)=2^s\sum_{r=s}^{2n+1}a_r\frac{r!}{(r-s)!}i^rt^r+o(t^{2n+1}),
\]
and thus
\begin{equation*}
\begin{split}
(\omega-\overline\omega)^{n+p-\ell}\overline{a_{p}(\omega)}b^{(n-\ell)}(\omega)=\frac{(-1)^p}{p!}2^{n+p-\ell}&\left(\sum_{r=n-\ell}^{2n+1}a_r\frac{r!}{(r-n+\ell)!}i^rt^r+o(t^{2n+1})\right)\\
&\times\left(\sum_{j=p}^{2n+1}\overline{a_j}\frac{j!}{(j-p)!}(-i)^jt^j+o(t^{2n+1})\right).
\end{split}
\end{equation*}
We deduce from (\ref{eq1:asymptot-norme}) that
\begin{equation*}
\begin{split}
\|k_{\omega,n}^b\|_b^2=\dfrac{(-1)^nn!}{2^{2n+2}\pi}t^{-2n-1}\biggl[&(-1)^n\dfrac{(2n)!}{n!}-n!\sum_{p=0}^n\sum_{\ell=0}^n(-1)^{p+\ell}2^{n+p-\ell}\binom{n-p+\ell}{\ell}\\
&\times\left(\sum_{r=n-\ell}^{2n+1}a_r\binom{r}{n-\ell}i^rt^r\right)\left(\sum_{j=p}^{2n+1}\overline{a_j}\binom{j}{p}(-i)^jt^j\right)+o(t^{2n+1})\biggr],
\end{split}
\end{equation*}
and denoting by $c_n=\dfrac{(-1)^nn!}{2^{2n+2}\pi}$, we can write
\begin{equation*}
\begin{split}
\|k_{\omega,n}^b\|_b^2=c_n t^{-2n-1}\biggl[&(-1)^n\dfrac{(2n)!}{n!}-n!\sum_{s=0}^{2n+1}\lambda_{s,n}t^s+o(t^{2n+1})\biggr],
\end{split}
\end{equation*}
with
\[
\lambda_{s,n}:=i^s\sum_{p=0}^n\sum_{\ell=0}^n(-1)^{p+\ell}2^{n+p-\ell}\binom{n-p+\ell}{\ell}\sum_{r=0}^{s}\binom{r}{n-\ell}\binom{s-r}{p}(-1)^{s-r}a_r\overline{a_{s-r}},
\]
where we assumed that $\binom{a}{b}=0$ if $a<b$ or $b<0$.

Now we recall that $k_{\omega,n}^b$ is weakly convergent as $\omega$ tends radially to $x_0$ and thus $\|k_{\omega,n}^b\|_b$ remains bounded. Therefore we necessarily have
\begin{equation}\label{eq2:derivees}
(-1)^n\dfrac{(2n)!}{n!}-n!\lambda_{0,n}=0\qquad\hbox{and}\qquad \lambda_{s,n}=0,\quad (1\leq s\leq 2n),
\end{equation}
which implies that
\begin{equation}\label{eq2:convergence-en-norme}
\|k_{\omega,n}^b\|_b^2=-n!c_n\lambda_{2n+1,n}+o(1),
\end{equation}
as $\omega$ tends radially to $x_0$. But
\begin{align*}
\lambda_{2n+1,n}&=(-1)^n i\sum_{r=0}^{2n+1}(-1)^{r+1}a_r\overline{a_{2n+1-r}}\sum_{p=0}^n\sum_{\ell=0}^n(-1)^{p+\ell}2^{n+p-\ell}\binom{n-p+\ell}{\ell}\binom{r}{n-\ell}\binom{2n+1-r}{p}\\
&=(-1)^ni2^n\sum_{r=0}^{2n+1}A_{n,r}a_r\overline{a_{2n+1-r}}.
\end{align*}
According to Proposition~\ref{Prop:Fred}, we have $A_{n,r}=-2^n$ if $0\leq r\leq n$ and $A_{n,r}=2^n$ if $n+1\leq r\leq 2n+1$. Then we obtain
\[
\lambda_{2n+1,n}=(-1)^ni2^{2n}\left(\sum_{r=n+1}^{2n+1}a_r\overline{a_{2n+1-r}}-\sum_{r=0}^{n}a_r\overline{a_{2n+1-r}}\right).
\]
Now note that
\[
\sum_{r=n+1}^{2n+1}a_r\overline{a_{2n+1-r}}=\sum_{r=0}^n\overline{a_r}a_{2n+1-r},
\]
which means
\[
\lambda_{2n+1,n}=(-1)^{n+1}2^{2n+1}\Im{\rm{m}}\left(\sum_{r=0}^{n}\overline{a_r}a_{2n+1-r}\right).
\]
But Proposition~\ref{Prop:cacul-de-norme} implies that
\[
\lambda_{2n+1,n}=(-1)^{n+1}2^{2n+2}\frac{\pi}{n!^2}\|k_{x_0,n}^b\|_b^2,
\]
and finally using (\ref{eq2:convergence-en-norme}) and the definition of $c_n$, we obtain
\[
\|k_{\omega,n}^b\|_b^2=\|k_{x_0,n}^b\|_b^2+o(1),
\]
which proves that $\|k_{\omega,n}^b\|_b\longrightarrow\|k_{x_0,n}^b\|_b$ as $\omega$ tends radially to $x_0$. Since $k_{\omega,n}^b$ tends also weakly to $k_{x_0,n}^b$ in $\HH(b)$ as $\omega$ tends radially to $x_0$, we get the desired conclusion.

\endpf

\begin{remark}
\rm{We have already seen in the proof of Theorem \ref{thm:main} that if $x_0$ satisfies the condition (\ref{eq:condition-type-ahern-clark}), then $|a_0|=1$ and
\begin{equation}\label{eq:derives}
\sum_{p=0}^ka_p\overline{a_{k-p}}=0,\qquad (1\leq k\leq n),
\end{equation}
where $a_p:=\frac{b^{(p)}(x_0)}{p!}$. In fact, we can prove that the
relation (\ref{eq:derives}) is also valid for $n+1\leq k\leq 2n$,
$k$ even. Indeed, according to (\ref{eq2:derivees}), for $n+1\leq
s\leq 2n$, we have
\begin{align}\label{eq:toto1}
0=\lambda_{s,n}:=&i^s\sum_{p=0}^n\sum_{\ell=0}^n(-1)^{p+\ell}2^{n+p-\ell}\binom{n-p+\ell}{\ell}\sum_{r=0}^{s}\binom{r}{n-\ell}\binom{s-r}{p}(-1)^{s-r}a_r\overline{a_{s-r}}\nonumber\\
&=(-i)^s2^n\sum_{r=0}^s(-1)^r a_r\overline{a_{s-r}}\sum_{p=0}^n\sum_{\ell=0}^n(-1)^{p+\ell}2^{p-\ell}\binom{n-p+\ell}{\ell}\binom{r}{n-\ell}\binom{s-r}{p}\nonumber\\
&=(-i)^s2^n\sum_{r=0}^sa_r\overline{a_{s-r}}A_{n,r,s},
\end{align}
with
\begin{equation}\label{eq:toto2}
A_{n,r,s}:=(-1)^{r}\sum_{p=0}^n\sum_{\ell=0}^n(-1)^{p+\ell}2^{p-\ell}\binom{n-p+\ell}{\ell}\binom{r}{n-\ell}\binom{s-r}{p}.
\end{equation}
Using similar arguments as in the proof of Proposition \ref{Prop:Fred}, we show that for every $0\leq r\leq n<s$, we have
\begin{equation}\label{eq:toto3}
A_{n,r,s}=\binom{s}{n}\dfrac{\Gamma\left(\tfrac{s-n+1}{2}\right)\Gamma(\tfrac{s-n+2}{2})}{\Gamma\left(\tfrac{s-2n+1}{2}\right)\Gamma\left(\tfrac{s+2}{2}\right)}.
\end{equation}
In the proof of this identity, we use Bayley's Theorem
\cite{Andrews} which says that
\[
{}_2F_1\!{\scriptsize{\left[\begin{matrix}a,1-a\\
b\end{matrix};\tfrac{1}{2}\right]}}=\dfrac{\Gamma\left(\tfrac{b}{2}\right)\Gamma(\tfrac{1+b}{2})}{\Gamma\left(\tfrac{a+b}{2}\right)\Gamma\left(\tfrac{1-a+b}{2}\right)}\,.
\]
Now using (\ref{eq:toto2}) it is easy to see that $A_{n,s-r,s}=(-1)^sA_{n,r,s}$, and with (\ref{eq:toto1}) and (\ref{eq:toto3}), we obtain
\[
 \binom{s}{n}\dfrac{\Gamma\left(\tfrac{s-n+1}{2}\right)\Gamma(\tfrac{s-n+2}{2})}{\Gamma\left(\tfrac{s-2n+1}{2}\right)\Gamma\left(\tfrac{s+2}{2}\right)}\left(
\sum_{r=0}^na_r\overline{a_{s-r}}+(-1)^s\sum_{r=n+1}^sa_r\overline{a_{s-r}}
\right)=0.
\]
Now recall that the Gamma function is a meromorphic function in the complex plane  without zeros and with poles at zero and the negative integers . Therefore we see that if $s$ is even ($n+1\leq s\leq 2n$), then
\[
\sum_{r=0}^s a_r\overline{a_{s-r}}=0.
\]
But if $s$ is odd ($n+1\leq s\leq 2n$), then $\frac{s-2n+1}2$ is
zero or a negative integer and, using this argument, we are not able
to conclude that $\displaystyle\sum_{r=0}^s
a_r\overline{a_{s-r}}=0$. This still remains as an open question.}
\end{remark}

{\bf Acknowledgments:} We would like to thank F. Jouhet, E. Mosaki and J. Zeng for helpful discussions concerning the proof of Proposition \ref{Prop:Fred} and A. Baranov for some suggestions to improve the redaction of this text. A part of this work was done while the first author was visiting McGill University. He would like to thank this institution for its warm hospitality.

\bibliographystyle{acm}
\nocite{*}
\bibliography{biblio}

\begin{thebibliography}{10}

\bibitem{ak70}
{\sc Ahern, P.~R., and Clark, D.~N.}
\newblock Radial limits and invariant subspaces.
\newblock {\em Amer. J. Math. 92\/} (1970), 332--342.

\bibitem{ak71}
{\sc Ahern, P.~R., and Clark, D.~N.}
\newblock Radial {$n{\rm th}$} derivatives of {B}laschke products.
\newblock {\em Math. Scand. 28\/} (1971), 189--201.

\bibitem{AR}
{\sc Anderson, J.~M., and Rovnyak, J.}
\newblock On generalized {S}chwarz-{P}ick estimates.
\newblock {\em Mathematika 53}, 1 (2006), 161--168 (2007).

\bibitem{Andrews}
{\sc Andrews, G.~E., Askey, R., and Roy, R.}
\newblock {\em Special functions}, vol.~71 of {\em Encyclopedia of Mathematics
  and its Applications}.
\newblock Cambridge University Press, Cambridge, 1999.

\bibitem{Baranov05}
{\sc Baranov, A.~D.}
\newblock Bernstein-type inequalities for shift-coinvariant subspaces and their
  applications to {C}arleson embeddings.
\newblock {\em J. Funct. Anal. 223}, 1 (2005), 116--146.

\bibitem{Bolotnikov}
{\sc Bolotnikov, V., and Kheifets, A.}
\newblock A higher order analogue of the {C}arath\'eodory-{J}ulia theorem.
\newblock {\em J. Funct. Anal. 237}, 1 (2006), 350--371.

\bibitem{de-branges1}
{\sc de~Branges, L., and Rovnyak, J.}
\newblock Canonical models in quantum scattering theory.
\newblock In {\em Perturbation Theory and its Applications in Quantum Mechanics
  (Proc. Adv. Sem. Math. Res. Center, U.S. Army, Theoret. Chem. Inst., Univ. of
  Wisconsin, Madison, Wis., 1965)}. Wiley, New York, 1966, pp.~295--392.

\bibitem{de-branges2}
{\sc de~Branges, L., and Rovnyak, J.}
\newblock {\em Square summable power series}.
\newblock Holt, Rinehart and Winston, New York, 1966.

\bibitem{pD70}
{\sc Duren, P.~L.}
\newblock {\em Theory of {$H\sp{p}$} spaces}.
\newblock Pure and Applied Mathematics, Vol. 38. Academic Press, New York,
  1970.

\bibitem{Dyak91}
{\sc Dyakonov, K.~M.}
\newblock Entire functions of exponential type and model subspaces in {$H\sp
  p$}.
\newblock {\em Zap. Nauchn. Sem. Leningrad. Otdel. Mat. Inst. Steklov. (LOMI)
  190}, Issled. po Linein. Oper. i Teor. Funktsii. 19 (1991), 81--100, 186.

\bibitem{Dyak02}
{\sc Dyakonov, K.~M.}
\newblock Differentiation in star-invariant subspaces. {I}. {B}oundedness and
  compactness.
\newblock {\em J. Funct. Anal. 192}, 2 (2002), 364--386.

\bibitem{fatou}
{\sc Fatou, P.}
\newblock S\'eries trigonom\'etriques et s\'eries de {T}aylor.
\newblock {\em Acta Math. 30}, 1 (1906), 335--400.

\bibitem{Fricain}
{\sc Fricain, E.}
\newblock Bases of reproducing kernels in de {B}ranges spaces.
\newblock {\em J. Funct. Anal. 226}, 2 (2005), 373--405.

\bibitem{Fricain-Mashreghi}
{\sc Fricain, E., and Mashreghi, J.}
\newblock Boundary behavior of functions of the de branges-rovnyak spaces.
\newblock {\em Complex Analysis and Operator Theory\/} (to appear).

\bibitem{HSS}
{\sc Hartmann, A., Sarason, D., and Seip, K.}
\newblock Surjective {T}oeplitz operators.
\newblock {\em Acta Sci. Math. (Szeged) 70}, 3-4 (2004), 609--621.

\bibitem{helson}
{\sc Helson, H.}
\newblock {\em Lectures on invariant subspaces}.
\newblock Academic Press, New York, 1964.

\bibitem{Jury}
{\sc Jury, M.~T.}
\newblock Reproducing kernels, de {B}ranges-{R}ovnyak spaces, and norms of
  weighted composition operators.
\newblock {\em Proc. Amer. Math. Soc. 135}, 11 (2007), 3669--3675 (electronic).

\bibitem{sarason}
{\sc Sarason, D.}
\newblock {\em Sub-{H}ardy {H}ilbert spaces in the unit disk}.
\newblock University of Arkansas Lecture Notes in the Mathematical Sciences,
  10. John Wiley \& Sons Inc., New York, 1994.
\newblock A Wiley-Interscience Publication.

\bibitem{Shapiro1}
{\sc Shapiro, J.~E.}
\newblock Relative angular derivatives.
\newblock {\em J. Operator Theory 46}, 2 (2001), 265--280.

\bibitem{Shapiro2}
{\sc Shapiro, J.~E.}
\newblock More relative angular derivatives.
\newblock {\em J. Operator Theory 49}, 1 (2003), 85--97.

\bibitem{Slater}
{\sc Slater, L.~J.}
\newblock {\em Generalized hypergeometric functions}.
\newblock Cambridge University Press, Cambridge, 1966.

\end{thebibliography}

\end{document}